\documentclass[a4paper,12pt]{article}

\usepackage{xcolor}
\usepackage{amsmath}
\usepackage{amsfonts}
\usepackage{amssymb}
\usepackage{amsthm}
\usepackage{amsopn}
\usepackage{amscd}
\usepackage{pstricks}
\usepackage[utf8]{inputenc}
\usepackage{courier}
\usepackage{url}
\usepackage{color}
\usepackage{subfigure}
\usepackage{xypic}
\usepackage{graphicx}
\usepackage{paralist}
\usepackage{subfigure}
\usepackage[vcentering,dvips]{geometry}
\def\address#1#2{\begingroup
\noindent\parbox[t]{7.8cm}{%
\small{\scshape\ignorespaces#1}\par\vskip1ex
\noindent\small{\itshape E-mail address}%
\/: #2\par\vskip4ex}\hfill%
\endgroup}%

\author{Nathan Owen Ilten \& Lars Kastner}
\title{Calculating Generators of Multigraded Algebras}
\date{}

\definecolor{lcolor}{rgb}{1.0,0.3,0.0}
\newcommand{\lars}[1]{}

\newcommand{\CC}{\mathbb{C}}
\newcommand{\QQ}{\mathbb{Q}}
\newcommand{\ZZ}{\mathbb{Z}}
\newcommand{\PP}{\mathbb{P}}
\newcommand{\G}{\mathcal G}
\newcommand{\D}{\mathcal D}
\newcommand{\A}{\mathcal A}
\newcommand{\CO}{\mathcal O}

\newcommand{\wt}{\omega}

\newcommand{\ray}{\rho}

\newcommand{\dfan}{\mathcal{S}}

\newcommand{\adm}{\A(\D,M)}

\newcommand{\admone}{\A(\D,M')}
\newcommand{\adwim}{\A(\D|_{\wt^i},M)}
\newcommand{\adsm}{\A(\D|_{\sigma},M)}
\newcommand{\adsmone}{\A(\D|_{\sigma},M')}

\newcommand{\qfield}{\CC(Y)(M)}
\newcommand{\quot}[2]{^{\displaystyle #1}\!\! \left/ \!\! _{\displaystyle #2}\right.}
\newcommand{\xvert}{\mathbf{vert}}
\newcommand{\xray}{\mathbf{ray}}
\newcommand{\tD}{\widetilde{\mathcal{D}}}

\newcommand{\tsigma}{\widetilde{\sigma}}

\DeclareMathOperator{\pos}{pos}
\DeclareMathOperator{\relint}{relint}
\DeclareMathOperator{\Spec}{Spec}

\DeclareMathOperator{\rk}{rk}

\DeclareMathOperator{\CDiv}{CaDiv}
\DeclareMathOperator{\Div}{div}
\DeclareMathOperator{\Quot}{Quot}
\DeclareMathOperator{\proj}{proj}
\DeclareMathOperator{\spec}{Spec}

\DeclareMathOperator{\conv}{conv}
\DeclareMathOperator{\Cl}{Cl}
\DeclareMathOperator{\Eff}{Eff}
\DeclareMathOperator{\Cox}{Cox}
\newcommand{\cox}{\textrm{cox}}

\theoremstyle{theorem}
\newtheorem{theorem}{Theorem}[section]

\theoremstyle{definition}

\newtheorem{alg}[theorem]{Algorithm}
\newtheorem{rem}[theorem]{Remark}
\newtheorem*{example}{Example}

\newcommand{\coeffb}{
\begin{pspicture}(-3,0)(3,3)
\pspolygon*[linecolor=lightgray](-2,2)(-1,1)(1,1)(2,2)
\psdots(-2,2)(-1,2)(0,2)(1,2)(2,2)(-1,1)(0,1)(1,1)
\psline{<->}(-2,2)(-1,1)(1,1)(2,2)
\rput(-1.5,.7){$(-1,1)$}
\rput(1.5,.7){$(1,1)$}
\end{pspicture}}

\newcommand{\coeffa}{
\begin{pspicture}(-3,0)(3,3)
\pspolygon*[linecolor=lightgray](-1.5,2)(0,.5)(1.5,2)
\psdots(-2,2)(-1,2)(0,2)(1,2)(2,2)(-1,1)(0,1)(1,1)
\psline{<->}(-1.5,2)(0,.5)(1.5,2)
\rput(0,.2){$(0,1/2)$}
\end{pspicture}}

\newcommand{\pdomain}{
\begin{pspicture}(-3,0)(3,3)
\psdots(-2,2)(-1,2)(0,2)(1,2)(2,2)(-1,1)(0,1)(1,1)(0,0)
\psline{<->}(-2,2)(0,0)(2,2)
\psline{->}(0,2)
\rput(-.5,1.5){$\wt^2$}
\rput(.5,1.5){$\wt^1$}
\end{pspicture}}

\newcommand{\deltacoeff}{
\begin{pspicture}(-3,0)(3,4)
\pspolygon*[linecolor=lightgray](-2.5,3)(-1,1.5)(2.5,1.5)(2.5,3)
\psdots(-1,2)(0,2)(1,2)(2,2)(-2,3)(-1,3)(0,3)(1,3)(2,3)
\psline{<->}(-2.5,3)(-1,1.5)(2.5,1.5)
\rput(-1,1){$(-1,3/2)$}
\end{pspicture}}

\begin{document}
\maketitle
\begin{abstract}
We present an algorithm to find generators of the multigraded algebra $\A$ associated to an arbitrary p-divisor $\D$ on some variety $Y$. A modified algorithm is also presented for the case where $Y$ admits a torus action. We demonstrate our algorithm by computing generators for the Cox ring of the smooth del Pezzo surface of degree $5$.
\end{abstract}

\section*{Introduction}
Let $Y$ be a semiprojective normal variety and $M$ some lattice with $M_\QQ$ the associated $\QQ$-vector space. Consider some full-dimensional polyhedral cone $\wt\subset M_\QQ$ with relative interior $\relint \wt$.
Then a \emph{p-divisor} on $Y$ with weight cone $\wt$ is a convex, piecewise linear function
$$
\D:\wt\to \CDiv_\QQ{Y}
$$
satisfying
\begin{enumerate}
	\item $\D(u)$ is \emph{semiample} for all $u\in \wt$, that is, a multiple of $\D(u)$ is globally generated;
	\item $\D(u)$ is \emph{big} for all $u \in \relint \wt$, that is, a multiple of $\D(u)$ admits a section with affine complement.
\end{enumerate}
To such a p-divisor $\D$ we can associate a multigraded section algebra:
$$
\adm:=\bigoplus_{u\in \wt\cap M} H^0\big(Y,\CO(\D(u))\big)\cdot\chi^u
$$
By a result of K. Altmann and J. Hausen, $\adm$ is an integral, normal, finitely generated $M$-graded $\CC$-algebra, see \cite[Theorem 3.1]{altmann:06a}. Furthermore, by Theorem 3.4 of the same article, all such algebras arise in this fashion.
This correspondence generalizes that between rational polyhedral cones and toric algebras; the geometry of $\Spec \adm$ can be readily studied in terms of the geometry of $Y$ together with the p-divisor $\D$.

Although K. Altmann and J. Hausen prove that $\adm$ is finitely generated, their proof gives no explicit method for finding a generating set.  The problem of computing generators has already been dealt with in some simple cases:
\begin{enumerate}
	\item $\dim Y=0$, i.e. $Y$ is a point. In this case, $\D$ is the trivial map, and $\adm$ is a normal toric algebra:
		$$
\adm=\bigoplus_{u\in \wt\cap M} \CC\cdot\chi^u=\CC[\wt\cap M]
		$$
		Generators of $\adm$ are given by calculating generators of the semigroup $\wt\cap M$. If $\wt$ is pointed, there is in fact a unique minimal set of semigroup generators called a \emph{Hilbert basis}, which can be calculated using a number of algorithms, see for example \cite{hemmecke:02a} or \cite{normaliz}.
	\item $\dim Y=1$, i.e. $Y$ is a smooth curve. This case was dealt with by H. S\"u\ss{} and the first author in \cite{ilten:09d}. The strategy here is to subdivide $\wt$ into cones on which $\D$ is linear, and then to use bounds on the degrees of the $\D(u)$ to find a finite set of weights $\G(\D)\subset \wt\cap M$ in which $\adm$ is generated. The set $\G(\D)$ can in fact be calculated without knowing anything about the geometry of $Y$ beyond its genus.
\end{enumerate}

We present an algorithm to compute generators of $\adm$ for an arbitrary p-divisor $\D$. In contrast to the simple cases above, in order to compute a set of weights in which $\adm$ is generated, we will need more specific geometric information about $Y$.  In particular, our algorithm requires the ability to calculate global sections of any divisor on $Y$.

We also present a modified algorithm for a special case. Indeed, suppose that the variety $Y$ admits a torus action. Then up to an integral closure, $\adm$ can be written as the union of certain subalgebras $\A_i$ such that the torus action on $Y$ lifts to one on $\Spec \A_i$. Instead of calculating global sections of divisors on $Y$, we can instead  calculate global sections of divisors on a suitable quotient $Z$ of $Y$. By imposing some integrality conditions on $\D$, we can reduce to calculating generators of $\A(\D_i,Z)$, where $\D_i$ are certain p-divisors on $Z$.

Our algorithms fit into the box of tools being developed to study $T$-varieties, see \cite{altmann:06a}, \cite{altmann:08a}, and \cite{tsurvey}. A major motivation for us in developing our algorithms was to find a method for effectively computing generators of Cox rings of log del Pezzo surfaces. Such rings are described in terms of a p-divisors in \cite{altmann:09a}. As an example, we show how our algorithm may be used to compute generators for the Cox ring of $S_5$, the smooth del Pezzo surface of degree $5$. In future work, we hope to compute generators of Cox rings of other log del Pezzo surfaces for which generators are not yet known.   

We now describe the organization of this article. In Section \ref{sec:prelim} we show how to reduce the general problem to some simpler cases, and  state a result by O. Zariski which is essential for our general algorithm. In Section \ref{sec:gen} we present and prove the effectiveness of our algorithm for a general p-divisor. Section \ref{sec:opt} contains a discussion on implementation and improvement possibilities. We show how our algorithm may be used to find the generators of the Cox ring of $S_5$ in Section \ref{sec:ex}. Finally, in Section \ref{sec:action} we present a modified algorithm for the case that $Y$ admits a torus action.

\section{Some Preliminaries}\label{sec:prelim}
Let $Y$ be semiprojective and normal and $A=H^0(\CO_Y,Y)$.
Consider any p-divisor $\D:\wt\to\CDiv_\QQ Y$.  For any subcone $\wt'\subset\wt$, we can consider $\D|_{\wt'}$, the restriction of $\D$ to $\wt'$. Then $\D|_{\wt'}$ is clearly also a p-divisor.  Furthermore, if $\{\wt^1,\ldots,\wt^r\}$ is some set of cones which covers $\wt$, then all homogeneous elements of $\adm$ are contained in at least one $\adwim$, so we have
$$
\adm=\sum_{i=1}^r \adwim.
$$
Thus, to find generators of $\adm$, it will suffice to find generators of $ \adwim$ for all $i$.
Since $\D$ is piecewise linear on $\wt$, we can in fact find such cones $\wt^i$ covering $\wt$ with $\D|_{\wt^i}$ linear. This means that after passing to some restriction of $\D$, we can actually reduce to the case that $\D$ is linear and $\wt$ is simplicial. 

Likewise, consider a sublattice $M'\subseteq M$ such that $\rk{M'}=\rk{M}$. We then have that
$$
\adm\subseteq\overline{\admone}^{\qfield}
$$
where the term on the right is the integral closure of $\admone$ in the function field $\qfield$.
Using $\admone\subseteq\adm$ and normality of $\adm$ in its quotient field $\qfield$ we even obtain equality. Thus, we can calculate generators of $\adm$ by first calculating generators for $\admone$ and then calculating generators for a certain integral closure.

\begin{rem}\label{rem:pdiv}
A p-divisor can always be represented as a formal sum over prime divisors on $Y$ with polyhedral coefficients
$$
\D=\sum_{P\subset Y} \D_P\otimes P
$$
where for each prime divisor $P$, $\D_P$ is a polytope in $M_\QQ^*$ with tail cone $\wt^\vee$, and $\D_P$ differs from $\wt^\vee$ for only finitely many $P$. To get a map $\D:\wt\to\CDiv_\QQ Y$, we set
$$
\D(u):=\sum_{P\in Y} \min \langle \D_P, u \rangle \cdot P
$$
for any $u\in\wt$.
As a convention, we may omit any prime divisor appearing in the representation of $\D$ whose coefficient is $\wt^\vee$.

Given such a representation of $\D$, we can find subcones $\wt'\subset\wt$ on which $\D$ is linear by considering elements of the coarsest common refinement of the normal fans of the $\D_P$.
\end{rem}

\begin{figure}
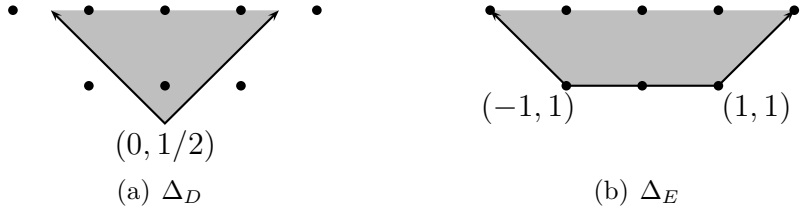
\begin{center}
	\subfigure[$\Delta_D$]{\coeffa}
	\subfigure[$\Delta_E$]{\coeffb}
\end{center}
	\caption{Polyhedral coefficients for a p-divisor}\label{fig:ex1}
\end{figure}

\begin{figure}
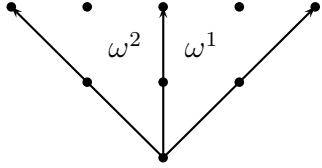

	\begin{center}
		\pdomain
	\end{center}
	\caption{Domain of a p-divisor}\label{fig:ex2}
\end{figure}

\begin{example}[A p-divisor on $\PP^2$] 
Let $M=\ZZ^2$, $Y=\PP^2=\proj [x,y,z]$, and consider the divisors $D=V(xyz)$ and  $E=V( (y-z)(x-z)(x-y))$. Let 
$$\D=\Delta_D \otimes D+\Delta_E\otimes E,$$
where $\Delta_D$  and $\Delta_E$ are the polytopes in $M^*_\QQ$ pictured in Figure \ref{fig:ex1}. Then $\D$ is a p-divisor. Its domain $\wt$ is pictured in Figure \ref{fig:ex2} together with the subcones $\wt^1,\wt^2$ on which it is linear. Explicitly, we have
\begin{align*}
\D:\wt&\to\CDiv_\QQ\PP^2\\
u=(u_1,u_2)&\mapsto\begin{cases}
	1/2u_2D+(u_2-u_1)E &\mbox{if }u\in\wt_1\\
	1/2u_2D+(u_2+u_1)E &\mbox{if }u\in\wt_2.\\
\end{cases}
\end{align*}
\end{example}	
The main ingredient we will need for our algorithm is the following result of Zariski, which tells us how to get generators of $\adm$ in a fairly simple case (assuming that we can compute a certain integral closure). Let $\D$ be a linear p-divisor on $\wt$, where $\wt=\QQ^n$. Assume furthermore that $\D(e_i)$ is an effective and globally generated $\ZZ$-divisor for $1\leq i \leq n$, where $e_1,\ldots,e_n$ is the standard basis of $\QQ^n$.
\begin{theorem}[{cf. \cite[Theorem 4.2]{zariski:62a}}] \label{theorem:zar2}
	Let $R$ be the $A$-algebra generated by $$H^0(Y,\CO(\D(e_i)))\cdot \chi^{e_i}$$ for $1\leq i \leq n$. Then the integral closure of $R$ in $\CC(Y)(M)$ is $\adm$.\footnote{Zariski's theorem as originally stated makes no assumptions about the bigness of $\D(u)$. The theorem is also only proved for projective $Y$, but the proof can be adapted without problem for semiprojective $Y$.}
\end{theorem}

\begin{rem}We can use this result to compute generators even in the case where $\D(e_i)$ isn't effective. Indeed, let $\D$ be linear on $\wt=\QQ^n$ as above with $\D(e_i)$ globally generated and integral but not necessarily effective. Then we can construct a new p-divisor $\D'$ satisfying $\D'(e_i)$ effective by setting
$$
\D'(u)=\D'(u_1,\ldots,u_n)=\sum u_i(\Div(s_i)+\D(e_i))
$$
for any global sections $s_i\in H^0(Y,\CO(\D(e_i)))$.
We then have
$$
\A(\D',M)\cong\adm
$$
with the isomorphism given by multiplying any homogeneous $s\cdot \chi^u\in \adm$ with $\prod s_i^{-u_i}$.  
\end{rem}

We finish this section with some notation and a final remark.
For any ray $\rho$ in $ M_\QQ$ or $M_\QQ^*$, we will identify $\rho$ with its primitive lattice generator. Likewise, for  any element $v$ of $M_\QQ$ or $M_\QQ^*$, let $\mu(v)$ be the smallest non-negative integer such that $\mu(v)v$ is a lattice point. 
\begin{rem}
	For any p-divisor $\D$, we have that $\Quot \adm=\qfield$. This follows from \cite[Theorem 3.1]{altmann:06a}, which can be used to show that $\spec \adm$ is birational to $Y\times \spec \CC[M]$. 
\end{rem}

\section{The Algorithm} \label{sec:gen}
The discussion in the previous section should already hint at our strategy for calculating generators of $\adm$. We will first reduce to the case that $\D$ is linear and integral, and then calculate the generators mentioned in Theorem \ref{theorem:zar2}. To actually get generators for $\adm$, we need to calculate an integral closure of the resulting algebra. By judiciously adding elements to this algebra, we reduce to the problem of calculating its normalization.
\begin{alg}\label{algo:gen}
\begin{itemize}{\bf (Computing generators in general case)}
\item[] {\bf INPUT:} $Y$ a normal semiprojective variety, $\D$ a p-divisor on $Y$.
\item[] {\bf OUTPUT:} $G\subseteq \CC(Y)[M]$ a set of generators of the algebra $\adm$.
\item[] {\bf PROCEDURE:}
\newcounter{lc}
\setcounter{lc}{0}
\begin{itemize}
\item[\refstepcounter{lc}\label{alg:poly}\arabic{lc}.] Find a polyhedral subdivision of $\wt$ with maximal cones $\{\wt^1,\ldots,\wt^l\}$ such that $\D|_{\wt^i}$ becomes linear for each $i$. 
\item[\refstepcounter{lc}\arabic{lc}.] Initialize an empty list $L$.
\item[\refstepcounter{lc}\arabic{lc}.] For each $i=1,\ldots, l$ do
\begin{itemize}
\item[\refstepcounter{lc}\arabic{lc}.] For each generating ray $\ray\in M$ of $\wt^i$ do:
\begin{itemize} 
\item[\refstepcounter{lc}\label{alg:mult}\arabic{lc}.] Find $k_{\ray}$ such that $\D(k_{\ray}\ray)$ is base point free and integral by taking generators of $H^0(Y,\D(k\ray))$ (as an $A$-module) and intersecting the corresponding divisors for increasing $k$. Keep the generators in the cache for the next step;
\item[\refstepcounter{lc}\label{alg:gen}\arabic{lc}.] Take generators $\{\eta_1,\ldots,\eta_{s_{\ray}}\}$ of $H^0(Y,\D(k_{\ray}\ray))$ and add the elements $\eta_j\chi^{k_{\ray}\ray}$ to $L$;
\end{itemize}
\end{itemize}  
\item[\refstepcounter{lc}\label{alg:lattice}\arabic{lc}.] Take a lattice basis $\{b_1,\ldots,b_n\}$ of $M$ that is contained in the interior of $\wt$.
For each $b_k$, do:
\begin{itemize}
\item[\refstepcounter{lc}\label{alg:kray}\arabic{lc}.] Initialize an empty list $G$ and a counter $j=1$. While the greatest common divisor of the elements in $G$ does not equal $1$ do: (For $G=\emptyset$ we define $\gcd(G):=0$.)
\begin{itemize}
\item[\refstepcounter{lc}\arabic{lc}.] Compute $H^0(Y,\D(j\cdot b_k))$;
\item[\refstepcounter{lc}\arabic{lc}.] If $H^0(Y,\D(j\cdot b_k))\not=0$ take a non-zero element $s\in H^0(Y,\D(j\cdot b_k))$ and add $s\cdot\chi^{j\cdot b_k}$ to $L$ and add $j$ to $G$;
\item[\refstepcounter{lc}\label{alg:kray:end}\arabic{lc}.] Increase $j$ by one;
\end{itemize}
\end{itemize}
\item[\refstepcounter{lc}\label{alg:nor1}\arabic{lc}.] Fix a ray $\ray\in M$ in the relative interior of $\wt$. Let $j:=0$. While $\CC(Y)\not\subseteq\Quot(A[L])$ add generators of $H^0(Y,\D(j\ray))$ to $L$ and increase $j$ by one;
\item[\refstepcounter{lc}\arabic{lc}.] Normalize $A[L]$ and return the generators of this normalization.
\end{itemize}
\end{itemize}
\end{alg}
\begin{proof}[Proof of correctness and termination]
We first claim that the algorithm will finish in finite time. Indeed, since we assume all cones to be polyhedral, we only have finitely many rays to examine. The corresponding divisors are semiample, thus steps \ref{alg:poly} through \ref{alg:gen} terminate in finite time.
The loops in step \ref{alg:lattice} and \ref{alg:nor1} terminate since for each ray in the interior of $\wt$ we get a p-divisor whose corresponding subalgebra of $\adm$ is finitely generated.
Now the normalization of a finitely generated algebra again is finitely generated, hence it can be computed in finite time. Thus, termination is guaranteed.

We now show that the algorithm yields the desired result.
Let $s\chi^u\in\adm$ with $s\in H^0(Y,\D(u))$ a homogeneous element of $\adm$. Now $u$ is contained in some $\wt^i$ of the polyhedral subdivision from step \ref{alg:poly}, hence $s\cdot\chi^u\in\adwim$. Fix $\rk{M}$ generating rays of $\wt^i$ such that $u$ is contained in the subcone $\sigma\subseteq\wt^i$ built by these rays. Then $s\chi^u\in\adsm$ holds. For each generating ray $\ray\in M$ of $\sigma$ find $k_{\rho}$ as described in step \ref{alg:kray} and let $M'$ the lattice generated by the $k_{\ray}\ray$. Then $s\chi^u$ is contained in the integral closure of $\adsmone$ in $\qfield$. Now we are in the setting of Theorem \ref{theorem:zar2} and can compute generators for a ring whose integral closure in $\CC(Y)(M')$ is exactly $\adsmone$. These generators are added to $L$ in step \ref{alg:gen}. Thus $s\chi^u$ is contained in the integral closure of $A[L]$ in $\qfield$.

Steps \ref{alg:kray} through \ref{alg:nor1} add elements to $L$ to make sure the quotient field of $A[L]$ becomes $\qfield$. Hence $s\chi^u$ is contained in the normalization of $A[L]$. This implies that $\adm$ is contained in the normalization of $A[L]$, since the former is generated by its homogeneous elements. The opposite inclusion is immediate, since every element of $A[L]$ is contained in $\adm$ and $\adm$ is normal.\par

\end{proof}

\begin{example}[A p-divisor on $\PP^2$]
	We continue the example from Section \ref{sec:prelim}. Figure \ref{fig:ex2} already pictures  a subdivision of $\wt$ on which $\D$ becomes linear. We collect the rays:
$$
\ray_1:=\left(\begin{array}{c} -1 \\ 1\end{array} \right),\quad
\ray_2:=\left(\begin{array}{c} 0 \\ 1\end{array} \right),\quad
\ray_3:=\left(\begin{array}{c} 1 \\ 1\end{array} \right).
$$
The divisors corresponding to each ray are:
$$
\D(\ray_1) = 1/2\cdot D,\quad
\D(\ray_2) = 1/2\cdot D + E,\quad
\D(\ray_3) = 1/2\cdot D.
$$
Now for step \ref{alg:mult} we want to determine $k_{\ray_i}$ for each ray $\ray_i$. On $Y=\PP^2$, every (integral) effective divisor automatically is semiample, so we only need to ensure integrality of the corresponding divisors.
Hence $k_{\ray_i}=2$ for $i=1,2,3$.
For $i\in\ZZ_{\ge 0}$ let $\CC[x,y,z]_i$ denote the $\CC$-vectorspace of homogeneous polynomials of degree $i$. We now compute the global sections of $\D(k_{\ray_i}\cdot\ray_i)$ as subvectorspaces of $\CC(Y)=\CC[x,y,z]{( (0) )}$, the homogeneous localization of $\CC[x,y,z]$ at $(0)$:
\[
H^0(Y,\D(2\cdot\ray_1))=
H^0(Y,\D(2\cdot\ray_3))=
H^0(Y,D)=\frac{1}{xyz}\CC[x,y,z]_3
\] 
and
\[
H^0(Y,\D(2\cdot\ray_2)) = H^0(Y,D+2E) = \frac{1}{xyz(y-z)^2(x-z)^2(z-y)^2}{\CC[x,y,z]_9}.
\]
For step \ref{alg:gen} we compute bases of the above vector spaces. Since they are rather large, we will not write them down. The first vector space is generated by $10$ elements, the second by $55$ elements.
Next we consider the loop in step \ref{alg:lattice}. A lattice basis can be given by $b_1=\rho_2=(0,1)$ and $b_2=\rho_3=(1,1)$. Since both $\D(b_i))$ are effective, we add the elements $\chi^{(0,1)}$ and $\chi^{(1,1)}$ to our set of generators and are done with step \ref{alg:lattice}.

Now we have collected $77$ generators, which is a rather large set. The ideal of relations has $>2000$ generators which makes computing the normalization impossible from a practical standpoint.
Thus, we will reduce the set of generators as follows: Let $f$ be a generator. 
Now if $f^n$ is algebraically dependent on the other generators for some $n\in\ZZ_{>0}$ we will remove this $f$.
Since we will compute the integral closure later on we are sure to keep $f$ in our algebra. However the quotient field of the algebra generated by our elements collected so far might change, so before normalizing, we need to make sure that $\CC(Y)\subset\Quot(A[L])$.
Eliminating generators we obtain the following smaller list of generators:
\begin{align*}
\frac{x^3}{f_1}\chi^{(-2,2)},
\frac{y^3}{f_1}\chi^{(-2,2)},
\frac{z^3}{f_1}\chi^{(-2,2)},
\frac{x^9}{f_1f_2^2}\chi^{(0,2)},
\frac{y^9}{f_1f_2^2}\chi^{(0,2)},
\frac{z^9}{f_1f_2^2}\chi^{(0,2)},\\
\frac{x^3}{f_1}\chi^{(2,2)},
\frac{y^3}{f_1}\chi^{(2,2)},
\frac{z^3}{f_1}\chi^{(2,2)},
\chi^{(0,1)},\chi^{(1,1)},
\frac{x^2y}{f_1}\chi^{(-2,2)},
\frac{xy^2}{f_1}\chi^{(-2,2)}
\end{align*}
where 
\[
f_1:=xyz \mbox{ and }  f_2:=(y-z)(x-z)(x-y).
\]
The last three generators were added to ensure that the quotient field is indeed $\CC(Y)(M)$.
The algebra $\adm$ is just the normalization of the subalgebra of $\CC(Y)(M)$ generated by these thirteen generators. For a concrete description, see the continuation of this example in Section \ref{sec:action}.
\end{example}

\section{Implementation and Optimization}\label{sec:opt}
Existing mathematical software supplies most of the ingredients necessary for an implementation of the algorithm described above. All polyhedral computations can easily be handled with the program Polymake \cite{polymake} or the Macaulay2 package Polyhedra \cite{polyhedra}. A central requirement of the algorithm is the ability to compute global sections of divisors on semiprojective varieties. This already can be done with existing code in some special situations using Macaulay2 \cite{macaulay}; it would be a useful addition to this software system to create a package which can do this in general. Finally, the last step of the algorithm requires the calculation of a normalization. Both Macaulay2 \cite{macaulay} and Singular \cite{singular} provide the necessary tools to do this. 

The list $L\subset\qfield$ constructed by the algorithm might be significantly larger than necessary.
In order to compute the normalization of $A[L]$ in a computer algebra system, one must first compute a presentation of $A[L]$; the larger $L$ is, the more difficult this will become. Thus, in order to optimize the algorithm, one may want to eliminate unnecessary elements of $L$ while running the algorithm. We outline two strategies for doing this below.

In steps \ref{alg:kray} to \ref{alg:kray:end} we want to ascertain that the lattice $M_1$ generated by the weights $u$ of the elements $s\cdot\chi^u$ collected in $L$ so far in fact equals $M$. To reduce the number of elements we have to add one may compute a lattice basis $H$ of $M_1$ in Hermite normal form. Now we reduce each new $u$ for every new $s\cdot\chi^u$ by $H$ to check whether it is already contained in $M_1$. If $u$ is already contained in $M_1$ we will not need $s\cdot\chi^u$ for our completion procedure. If $u$ is not contained in $M_1$ we add the element $s\cdot\chi^u$ to $L$. Now the lattice $M_1$ changes and possibly becomes finer in which case we recompute $H$.

It is also possible to check whether elements found in step \ref{alg:gen} and \ref{alg:nor1} are already contained in the algebra generated by the previous elements. Depending on the number of elements already in the collection this might involve large computations of standard bases. Several algorithms for deciding algebraic dependence can be found in \cite{sturmfels:08}.

\section{Example: The Cox Ring of a del Pezzo Surface}\label{sec:ex}
	Let $Y$ be any $\QQ$-factorial projective variety with $\Cl(Y)$ finitely generated and free abelian. Then the Cox ring of $Y$ is 
	$$\Cox (Y)=\bigoplus_{ D\in\Cl(Y)}H^0(Y,D)$$
	with multiplication defined by a choice of basis of $\Cl(Y)$. If $\Cox (Y)$ is finitely generated, $Y$ is colloquially called a Mori Dream Space, and K. Altmann and J. Wisniewski have described a p-divisor $\D^\cox$ such that $\A(\D^\cox,\Cl(Y))=\Cox(Y)$  \cite{altmann:09a}. If $Y$ is a surface, then $\D^\cox$ in fact lives on $Y$. In the following, we use this description of $\D^\cox$ in conjunction with our algorithm to find generators of $\Cox(Y)$ for a special Mori Dream Surface $Y$. 

Let $Y=S_5$ be the smooth del Pezzo surface of degree $5$, i.e. the blow-up of $\PP^2$ at $4$ general points $P_1,\ldots,P_4$, for example $[1:0:0]$, $[0:1:0]$, $[0:0:1]$, $[1:1:1]$.
Let $M$ denote the following $(5\times 10)$-matrix:
\[
\left(
\begin{array}{cccccccccc}
0& 0& 0& 0& 1& 1& 1& 1& 1& 1\\
1& 0& 0& 0& -1& 0& 0& 0& -1& -1\\ 
0& 1& 0& 0& 0& 0& -1& -1& 0& -1\\ 
0& 0& 1& 0& 0& -1& 0& -1& -1& 0\\ 
0& 0& 0& 1& -1& -1& -1& 0& 0& 0
\end{array}
\right).
\]
Let $\wt\subseteq \QQ^5$ be the positive hull of the columns of $M$. The class group $\Cl(Y)$ of $Y$ is generated by $[H]$ and $[E_i]$, $i=1,\ldots,4$, $H$ being the pullback of a line of $\PP^2$ and $E_i$ the exceptional divisors. For $u\in\wt$ we define $d(u):=u_0\cdot H+\sum_{i=1}^4 u_i\cdot E_i$.
This yields a bijection of $\wt\cap\ZZ^5$ and $\Eff(Y)$ the cone of effective divisors on $Y$.
Let $E_{ij}$ denote the pullback of the line in $\PP^2$ through the points $P_i$ and $P_j$. We use the description given in \cite{altmann:09a} to construct the p-divisor $\D^{Cox}$ on $Y$ associated to the Cox ring of $Y$:
\[
\begin{array}{cccc}
\D^{Cox}: & \wt & \to & \CDiv(S)\\
& u & \mapsto & \D(u)
\end{array}
\]
where
\[
\D(u) = H.d(u)\otimes H + \sum_{i=1}^4\min(0,-E_i.d(u))\otimes E_i + \sum_{ij}\min(0,E_{ij}.d(u))\otimes E_{ij}.
\]
The Cox ring of $S_5$ can now be written as
\[
\Cox(S_5) = \A(\D^{Cox},\ZZ^5) = \bigoplus_{u\in\wt\cap\ZZ^5}H^0(S_5,\D(u))\cdot\chi^u \subseteq \CC(S_5)[\ZZ^5].
\]

To each of the appearing minima we associate the hyperplane $\{u\in\ZZ^5\ |\ \langle E, d(u)\rangle =0\}$. These are exactly the hyperplanes where the minima switch their values. Subdividing $\wt$ by intersecting with these hyperplanes we obtain a subdivision of $\wt$ into $241$ subcones on which $\D^{Cox}$ becomes linear. There are $160$ rays generating these cones. Thus, we have finished step one and at the same time found the rays needed for step three of the algorithm. 

Evaluating $\D^{Cox}$ on the primitive generators of these rays we obtain eleven different divisors up to linear equivalence, namely the zero divisor, $H$ , $2\cdot H$ , $H-E_i$ and $2\cdot H-2\cdot E_i$. As described in step five we check at what multiple each divisor becomes base point free. In fact, each of these divisors is itself already base point free.

Since the set of rays is fairly big and hence the set of generators will become even bigger, we want to make a small interruption and try to reduce the set of rays to consider. 
Assume we have $u^0,u^1,u^2\in\wt\cap\ZZ^5$ such that $u^0+u^1=u^2$, $\D(u^0) \sim_{lin} 0$ and $\D(u^1)\sim_{lin}\D(u^2)$. This means we automatically get surjectivity of the multiplication map 
\[
H^0(Y,\D(u^0))\times H^0(Y,\D(u^1))\to H^0(Y,\D(u^2)). 
\]
In terms of the algorithm this means that it is enough to collect the global sections for the rays $u^0$ and $u^1$ since every element of $H^0(Y,\D(u^2))\cdot\chi^{u^2}$ can be written as a product under the above multiplication map.
We use this statement to reduce the set of rays. \par
\vspace{0.5cm}
This leaves us with $23$ rays to consider for step six of the algorithm. We will now proceed by computing the global sections for each ray. Let
\[
P :=\ \quot{\CC[x_0,x_1,x_2,h,t_0,\ldots,t_9]}{(h\cdot f -1+\mbox{ toric relations})},
\]
where we associate to $t_i$ the $i$-th column of $M$ which also contains the Hilbert basis of $\wt$ and $f = x_0-x_1+x_2$ is the homogeneous equation in $\CC[x_0,x_1,x_2]$ defining $H$. The ideal of toric relations is generated by the binomials $t^u-t^v$ with $u,v\in\ZZ^{10}_{\ge 0}$ and $M\cdot (u-v)=0$. By the above convention we can write every $\chi^u$, $u\in\wt\cap\ZZ^5$ as a product of the $t_i$. Dividing out the toric relations makes the way of writing $\chi^u$ as a product unique.

Next we can write the global sections of the divisors associated to the remaining $23$ rays as elements of $P$. After eliminating algebraic dependencies from the resulting set of $57$ elements of $P$ we get the following $10$ generators of a subalgebra of $P$:
\begin{equation}\label{eq:gen}
	\begin{array}{c}
t_0,\ t_1,\ t_2,\ t_3,\ (x_1h-x_2h)\cdot t_4,\ (x_0h-x_1h)\cdot t_5,\\
(x_0h-x_2h)\cdot t_6,\ x_0ht_7,\ x_1ht_8,\ x_2ht_9.
\end{array}
\end{equation}
Thus, we have finished step six for all rays. Since every degree $t_i$ appears in the generating set,
we immediately finish with steps seven to eleven without adding any additional generators. 
Likewise, forming elements in the quotient field of the algebra of degree $0$ yields that $\CC(Y)$ is already contained in the quotient field, so we are immediately done with step twelve as well.
Finally, the subalgebra generated by this set of elements of $\CC(Y)[M]$ is already normal, so we may omit step thirteen.

Thus, we have shown that $\Cox(S_5)$ is generated by the elements listed in \eqref{eq:gen}.
One also observes that the coefficients of the $t_i$ appearing are exactly the $3\times3$-minors of the matrix
\[
\left(
\begin{array}{ccccc}
1 & 0 & 0 & 1 & x_0h\\
0 & 1 & 0 & 1 & x_1h\\
0 & 0 & 1 & 1 & x_2h
\end{array}
\right).
\]
Hence, we have obtained the generators of the Cox ring of $S_5$ as described in \cite{bp:04}.

Although a presentation of $\Cox(S_5)$ was in fact already known, there are numerous Mori Dream Surfaces for which the generators of the Cox rings are not yet known. We hope to perform similar computations for these examples in future work.

\section{Utilizing a Torus Action}\label{sec:action}
Let $Y$ now be a normal semiprojective variety admitting an effective action by an algebraic torus $T=\spec \CC[M']$. Then $Y$ may be reconstructed from a \emph{divisorial fan} $\dfan$ on a normal semiprojective $Z$: such a divisorial fan is a finite set of p-divisors on open subvarieties of $Z$ with respect to the lattice $M'$ satisfying some gluing conditions. The varieties $\{\spec \A(\D',M')\}_{\D'\in\dfan}$ glue together to give an affine open covering of $Y$, see \cite{altmann:08a} for more details. This generalizes the construction of arbitrary toric varieties via polyhedral fans.

Suppose that $\dfan$ is a divisorial fan on $Z$ which describes $Y$. One property of the p-divisors in $\dfan$ is that for any prime divisor $P\subset Z$, the set of polyhedral coefficients $\dfan_P=\{\D_P'\}_{\D'\in \dfan}$ form a polyhedral subdivision in $(M')^*_\QQ$.
In \cite{petersen:08a}, L. Petersen and H. S\"u\ss{} describe 
all $T$-invariant prime Weil divisors on $Y$ in terms of the vertices and rays of these polyhedral subdivisions. 
Indeed, there are ``vertical'' invariant prime divisors arising as the closure of a family 
of $\dim T$-dimensional $T$-orbits. 
Such divisors are parametrized by prime divisors 
$P\subset Y$ together with $v$ a vertex of $\dfan_P$ satisfying some 
additional condition \cite[Proposition 3.13]{petersen:08a}; let $\xvert_P(\dfan)$ denote the set of such vertices. We denote the divisor corresponding to such 
$P$ and $v$ by $D_{P,v}$. All other invariant prime divisors are ``horizontal'' and arise
as the closure of a family of $(\dim T-1)$-dimensional $T$-orbits. These are parametrized by rays $\rho$ of the tailfan of any of the $\dfan_P$ satisfying an
additional condition; let $\xray(\dfan)$ denote the set of such rays. 
We denote the divisor corresponding to a ray $r$ by $D_r$.

\begin{rem}
Often, we will want $Y$ to correspond to a divisorial fan $\dfan$ which is \emph{contraction free}. This means that any p-divisor $\D'\in\dfan$ lives on a smooth affine open subset of some  $Z$. This is not a serious restriction, since if we are given a p-divisor $\D$ on a $Y$ which is not of this form, we can blow up $Y$ and pull back $\D$ to arrive in this situation, see \cite[Remark 2.3]{upgrade}.
\end{rem}

Let $Y$ now be a normal semiprojective variety with an effective action by an algebraic torus $T$. Consider some $D\in\CDiv_\QQ Y$. We say that $D$ is \emph{$T$-moveable} if there exists $s\in H^0(Y,D)$ with $D+\Div(s)$ $T$ invariant. Let now $\D$ be a p-divisor $\D$ on $Y$ as in the introduction. We say that $\D$ is \emph{locally $T$-moveable} if $\D(u)$ is $T$-moveable for every $u\in\wt\cap M$. In what follows, we will see that if $\D$ is locally $T$-moveable, then we may calculate generators of $\adm$ by calculating generators for the algebras corresponding to certain p-divisors on quotients of $Y$ by $T$. In particular, if $\dim Y=\dim  T$, then the problem becomes completely combinatorial.

\begin{rem}
	Let $\D=\sum_{i=1}^k\Delta_i\otimes P_i$ be a representation of $\D$ as in Remark \ref{rem:pdiv}. If $Y$ admits an effective $T$-action and all $\Delta_i$ are lattice polyhedra, then $\D$ is automatically locally $T$-moveable.
\end{rem}

\begin{alg}
\begin{itemize}{\bf (Computing generators by using a $T$-action)}
\item[] {\bf INPUT:} $Y$ a normal semiprojective variety with effective $T$-action, $\D$ a p-divisor on $Y$ which is locally $T$-moveable.
\item[] {\bf OUTPUT:} $G\subseteq \CC(Y)[M]$ a set of generators of the algebra $\adm$.
\item[] {\bf PROCEDURE:}
\setcounter{lc}{0}
\begin{itemize}
	\item[\refstepcounter{lc}\label{step:subdivide}\arabic{lc}.] Find a polyhedral subdivision of $\wt$ with maximal cones $\{\wt^1,\ldots,\wt^l\}$ such that $\D|_{\wt^i}$ becomes linear for each $i$ and the $\wt^i$ are unimodular simplices. 
\item[\refstepcounter{lc}\arabic{lc}.] Initialize an empty list $L$.
\item[\refstepcounter{lc}\arabic{lc}.] For each $i=1,\ldots, l$ do

\begin{itemize}
	\item[\refstepcounter{lc}\label{step:sections}\arabic{lc}.] Find  $s_\rho\in H^0(Y,\D(\rho))$ for every ray $\rho$ of $\wt^i$ such that $\D(\rho)+\Div(s_{\rho})$ is $T$-invariant.
\item[\refstepcounter{lc}\arabic{lc}.] Blowup $Y$ to some $Y'$ such that with respect to the above $T$-action, $Y'$ can be represented by a contraction free divisorial fan $\dfan$ on some $Z$.
\item[\refstepcounter{lc}\arabic{lc}.]  Pullback $\D|_{\wt^i}$ to a p-divisor $\D'$ on $Y'$ and represent
	the p-divisor $\D''$
	$$
	u\mapsto \D'(u)+\sum_\rho u_\rho \Div(s_\rho)
	$$
	as
$$
\sum_{r\in\xray(\dfan)} \Delta_r \otimes D_r + \sum_{\substack{P\subset Y'\\v\in\xvert_P(\dfan)}} \Delta_{P,v} \otimes {\mu(v)}D_{P,v},
$$

	where in the first line, $\rho$ is ranging over all rays of $\wt^i$, and $u_\rho$ is the component of $u$ in the basis given by the rays $\rho$.
\item[\refstepcounter{lc}\arabic{lc}.]  Define a p-divisor $\tD = \sum \Delta_P \otimes P$ on $Z$ by
$$
\Delta_P = \conv\big\{\Delta_{P,v} \times \{v\} \mid v \in \xvert_P(\dfan) \big\} + \tsigma.
$$
and
$$
\tsigma = \pos \big\{ ( (\wt^i)^\vee \times \{0\}) \cup 
                      \bigcup_{r\in\xray(\dfan)} (\Delta_r \times \{r\}) \big\}.
$$

\item[\refstepcounter{lc}\label{step:lowerdim} \arabic{lc}.]Calculate  $M$-homogeneous generators $L_i\subseteq \CC(Y)[M] $ of $\A(\tD,M\times M')$. 
\item[\refstepcounter{lc}\arabic{lc}.] For each $f\chi^u\in L_i$, add $f\chi^u\prod_\rho s_\rho^{u_\rho}$ to $L$.
\end{itemize}
\item[\refstepcounter{lc}\arabic{lc}.] Return $L$. 
\end{itemize}
\end{itemize}
\end{alg}

\begin{proof}[Proof of correctness and termination]
	The proposed algorithm will terminate in finite time, assuming that sections demonstrating the local $T$-movability of $\D$ can be found effectively. Indeed, this assumption guarantees the termination of step \ref{step:sections}. Step \ref{step:lowerdim} will terminate by Algorithm \ref{algo:gen},  since $\tD$ truly is a p-divisor, see below. All other steps are finite polyhedral computations.

	Now fix some $i$ and consider the p-divisor $\D|_{\wt^i}$. Then $\A(\D,M)=\A(\D',M)$, and 
	$\A(\D'',M)$ is isomorphic to $\A(\D',M)$ via the map which sends $f\chi^u$ to $f\chi^u\prod_\rho s_\rho^{u_\rho}$.  Furthermore, by \cite[Theorem 2.2]{upgrade}, $\tD$ is a p-divisor and $\A(\D'',M)=\A(\tD,M\times M')$. Since $\A(\D,M)$ is generated by the subalgebras $\A(\D|_{\wt^i},M)$, this shows that the algorithm yields the correct result.
\end{proof}
\begin{rem}
	The above algorithm will function with a weaker assumption than $T$-moveability: indeed, we only need to assume that we can find a subdivision of $\wt$ as in step \ref{step:subdivide} such that for each $\wt^i$, we can (effectively) find an effective $T$-action on $Y$ together with sections $s_\rho$ as in step \ref{step:sections}.
\end{rem}
\begin{figure}
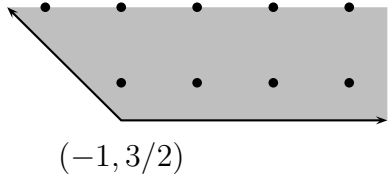

\begin{center}
\deltacoeff
\caption{The polyhedral coefficient $\Delta$}\label{fig:ex3}
\end{center}
\end{figure}

\begin{example}[A p-divisor on $\PP^2$]
We continue the example from Section \ref{sec:prelim}.
$Y=\PP^2$ is endowed with the standard action by $T=(\CC^*)^2$. 
Under this action, $Y$ may be described by the complete fan $\dfan$ in $\ZZ_\QQ^2$ with rays through $(1,0)$, $(0,1)$, and $(-1,-1)$.
The divisor $D$ is invariant under this action, but the divisor $E$ is not. Nonetheless, the p-divisor $\D$ is locally $T$-moveable.

For the first step of the algorithm, we may take the subdivision of $\wt$ into $\wt^1$ and $\wt^2$ pictured in Figure \ref{fig:ex2}. Let us focus on $\wt^1$, which has rays through $\rho_1=(0,1)$ and $\rho_2=(1,1)$. $\D(\rho_1)=1/2D+E$ and $\D(\rho_2)=1/2D$, so in step \ref{step:sections} we may take sections 
$$s_{\rho_1}=\frac{xyz}{(x-y)(x-z)(y-z)},\qquad s_{\rho_2}=1.$$

Since $\dfan$ lives just on a point, it is already contraction free. The p-divisor $\D''$ is just $(u_1,u_2)\mapsto (3/2u_2-u_1)D$. Denoting the three rays of $\dfan$ by $r_1,r_2,r_3$, this p-divisor may be represented as $\D''=\Delta\otimes(D_{r_1}+D_{r_2}+D_{r_3})$, where $\Delta$ is pictured in figure \ref{fig:ex3}.
The p-divisor $\tD$ is then just represented by the cone $\tsigma$ generated by the columns of the following matrix:
\begin{equation*}
\left(\begin{array}{c c c c c}
-1&1&-2&-2&-2\\
1&0&3&3&3\\
0&0&2&0&-2\\
0&0&0&2&-2
\end{array}\right)
\end{equation*}
To calculate homogeneous generators of $\A(\tD,M\times M')$, we calculate the Hilbert basis of the dual of $\tsigma$, which consists of the 66 columns of the following three matrices:
{\scriptsize
\begin{align*}
\left(\begin{array}{ccccccccccccccccccccccccccccccccccccccccc}
 0& 2& 1& 0& 2& 1& 0& 2& 1& 0& 0& 0& 2& 1& 0& 2& 1& 0& 0& 0& 0\\
 1& 2& 2& 2& 2& 2& 2& 2& 2& 2& 1& 1& 2& 2& 2& 2& 2& 2& 1& 1& 1\\ 
 -1& -1& -2& -3& 0& -1& -2& -1& -2& -3& 0& -1& 1& 0& -1& -1& -2& -3& 1& 0& -1\\ 
 -1& -1& -2& -3& -1& -2& -3& 0& -1& -2& -1& 0& -1& -2& -3& 1& 0& -1& -1& 0& 1\\
\end{array}\right)\\
\left(\begin{array}{ccccccccccccccccccccccccccccccccccccccccc}
2& 2& 2& 1& 0& 2& 1& 0& 0& 0& 0& 0& 1& 0& 1& 0& 1& 0& 1& 0& 1& 1\\ 
2& 2& 2& 2& 2& 2& 2& 2& 1& 1& 1& 1& 2& 2& 2& 2& 2& 2& 2& 2& 2& 2\\ 
1& 0& 2& 1& 0& -1& -2& -3& 2& 1& 0& -1& 2& 1& -2& -3& 3& 2& -2& -3& 3& 2\\
 0& 1& -1& -2& -3& 2& 1& 0& -1& 0& 1& 2& -2& -3& 2& 1& -2& -3& 3& 2& -1& 0\\
\end{array}\right)\\
\left(\begin{array}{ccccccccccccccccccccccccccccccccccccccccc}
1&1& 1& 1& 1& 0& 1& 0& 0& 0& 0& 0& 0& 0& 0& 0& 0& 0& 0& 0& 0& 0\\ 
1&2& 2& 2& 2& 2& 2& 2& 2& 2& 2& 2& 2& 2& 2& 2& 2& 2& 2& 2& 2& 2\\
0& 1& 0& -1& 4& 3& -2& -3& 4& -3& 5& -3& 5& 4& 3& 2& 1& 0& -1& -2& 6& -3\\ 
0&1& 2& 3& -2& -3& 4& 3& -3& 4& -3& 5& -2& -1& 0& 1& 2& 3& 4& 5& -3& 6\\
\end{array}\right)
\end{align*}}
To convert these into generators of $\A(\D'',M)$, we send a column $w=(w_1,w_2,w_3,w_4)$ to  
$$
(x/z)^{w_3}(y/z)^{w_4}\chi^{(w_1,w_2)}.
$$
Finally, to get generators of $\A(\D|_{\wt^1},M)$, we make the coordinate change
$$
\chi^{(0,1)}\mapsto\frac {xyz}{(x-y)(x-z)(y-z)}\chi^{(0,1)},\qquad \chi^{(1,1)}\mapsto\chi^{(1,1)}.
$$ 

Generators of $\A(\D|_{\wt^2},M)$ may be found similarly, showing that $\A(\D,M)$ can be generated by 132 generators found in the degrees $(0,1)$, $(0,2)$, $(\pm 1,1)$, $(\pm 1,2)$, $(\pm 2,2)$.
\end{example}
\begin{rem}
Although the p-divisor we started with in the above example appears fairly simple, describing the corresponding multigraded algebra in terms of a presentation is not. This underscores that fact that in many situations, a p-divisor is a much compacter way of encoding a multigraded algebra than via generators and relations.
\end{rem}

\bibliography{tgen}

\address{
Nathan Owen Ilten\\
Department of Mathematics\\
University of California\\
Berkeley, CA 94720\\
USA}{nilten@math.berkeley.edu}

\address{
Lars Kastner\\
Mathematisches Institut\\
Freie Universit\"at Berlin\\
Arnimallee 3\\
14195 Berlin, Germany}{kastner@math.fu-berlin.de}

\end{document}